# A Practical Method for Preventing Forced Wins in Ultimate Tic-Tac-Toe

Justin Diamond
11 July 2022

**Abstract.** Ultimate Tic-Tac-Toe is a variant of the popular Tic-Tac-Toe game. Two players compete to win three aligned "fields," with each field constituting its own miniature tic-tac-toe game. Each move determines which field the next player must play in. Prior studies have shown that there exists a forced winning strategy for the first player, whereby they can win in at least 29 moves and at most 43 moves.[1] This paper proposes a practical solution to the forced-win problem discovered by Bertholon et al., by putting forth a simple method for randomizing the first set of moves that are played. This method uses 5 randomly generated digits between 0 and 8 to arbitrarily place the first 4 moves of the game, and helps players avoid forced wins without having to change other rules of the game. This paper also investigates the probability that a random placement of the first 4 moves will lead to an easily calculable forced win, and shows that this probability is precisely 64/59049, or 0.108%.

**Introduction.** I will use the same set of rules and method of indexing as the study by Bertholon et al., such that all inquiry into Ultimate Tic-Tac-Toe (U3T) remains consistent. Please see their paper for an exhaustive explanation of U3T and a list of U3T rules.[2] Essential to this study is the fact that U3T is played on a board containing 9 fields arranged in a 3 × 3 grid. Since each field itself is further divided into 3 × 3 = 9 spots, the board contains 81 spots uniquely identifiable as $(i, j)$, where $i$ is the field and $j$ the spot. This indexing method is outlined in Figure 1 of the Bertholon et al. study (pictured below).[3]

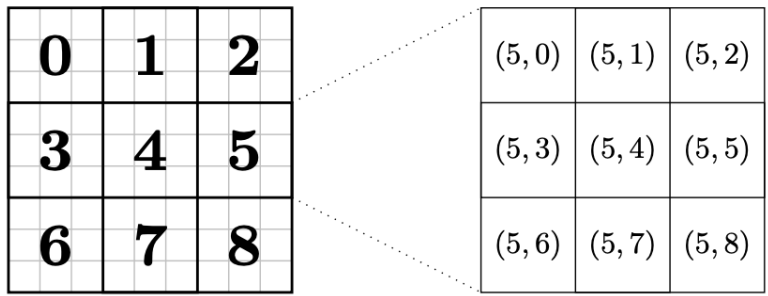

**Fig. 1.** Indexing for fields and spots. The board has 9 fields (*left*). Each field contains 9 spots, and we write $(i, j)$ for the $j$-th spot of the $i$-th field (*right*).

Because the board is divided in this way, a randomly generated digit between 0 and 8 can refer to any field, or spot within a field, that is present on the board. This makes it very easy to calculate and quantify the usefulness of the method for randomly generating starting positions which is proposed in this paper.

---

**Proposed Method.** A practical solution to the forced-win problem is randomizing the first set of moves that are played. By randomly generating 5 digits between 0 and 8, players can arbitrarily place the first 4 moves of the game. This helps players avoid forced wins without having to change other rules of the game. There is a slim probability that a random placement of the first 4 moves will result in an easily calculable forced win, or an illegal move.

The first randomly generated digit will refer to a field, $i$, and the remaining four digits will refer to a spot within a field, $j$. Here is how a randomly generated number such as 61245 would results in the first four placements, using the $(i, j)$ indexing method:

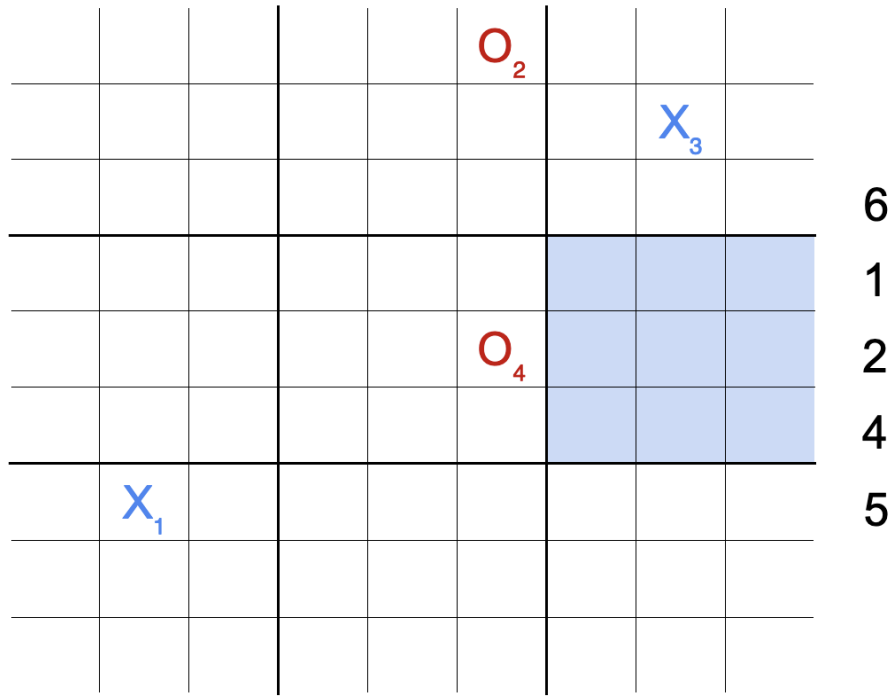

**Fig. 2.** An example of randomly generated placements for the 5-digit number 61245

1. 6 refers to the field that $X_1$ will be placed in.
2. 1 refers to the spot that $X_1$ will be placed in, and the field that $O_2$ will be placed in.
3. 2 refers to the spot that $O_2$ will be placed in, and the field that $X_3$ will be placed in.
4. 4 refers to the spot that $X_3$ will be placed in, and the field that $O_4$ will be placed in.
5. 5 refers to the spot that $O_4$ will be placed in, and the field that $X_5$, the first player-controlled move of the game will be placed in.

The resulting placements of 61245 are therefore:
    $X_1$: (6, 1)
    $O_2$: (1, 2)
    $X_3$: (2, 4)
    $O_4$: (4, 5)
    $X_5$: (5, $j$)

If the first 2 digits, and the 4th digit are all 4's, and the 3rd and 5th digit are non 4's, then an easily calculable forced win has been generated, per the strategy referenced by Bertholon et al..[4] An example of this could be the number 44148. The probability that this will occur is $(1/9) * (1/9) * (8/9) * (1/9) * (7/9)$ = 56/59,049, or a 0.0948% chance. If this scenario occurs, players can easily generate a new set of numbers.

It is also possible that an illegal move could be suggested by this method. For example, if the first 3 digits generated are 4, then an illegal move is suggested, since the center spot of the central field is already occupied by a move that has already been made. Other examples of illegal combinations are 84441, 81444, 44144, 41841, 84141. These scenarios are much more common than calculable forced wins, but are still outliers which can easily be avoided by generating a new number when they happen to occur.

**Conclusion.** This paper shows that the proposed method of randomly generating starting positions for Ultimate Tic-Tac-Toe is both practically easy and functionally useful. When 5 digits between 0 and 8 are randomly generated in order to derive a starting position, the odds of generating a position which can be exploited using strategies from Bertholon et al. are shown to be adequately small.[5] Although illegal starting positions or forced win positions may occasionally be suggested, the ease of generating new positions with a simple random number generator makes this method very appealing for players of U3T who wish to avoid theoretical forced wins.

---

[4] Bertholon, Guillaume; Géraud-Stewart, Rémi; Kugelmann, Axel; Lenoir, Théo; Naccache, David (June 3, 2020). "At Most 43 Moves, At Least 29: Optimal Strategies and Bounds for Ultimate Tic-Tac-Toe". arXiv:2006.02353v2
[5] Ibid.